\documentclass[draft]{amsart}
\usepackage{amsfonts}
\usepackage{amsmath}
\usepackage{amscd}
\usepackage{amssymb}
\usepackage{latexsym}

\theoremstyle{plain}
\newtheorem{stelling}{Theorem}

\swapnumbers
\newtheorem{theorem}[subsection]{Theorem}
\newtheorem{corollary}[subsection]{Corollary}
\newtheorem{lemma}[subsection]{Lemma}
\newtheorem{proposition}[subsection]{Proposition}

\newtheorem*{citedtheorem}{Theorem}

\theoremstyle{definition}
\newtheorem{definition}[subsection]{Definition}

\theoremstyle{remark}

\newtheorem{remark}[subsection]{Remark}

\swapnumbers

\newcommand\preprintnote {preprint on \myhomepage}
\newcommand\myhomepage{http://www.math.ohio-state.edu/\-\~{}schoutens}

\newcommand\OSU{\address{Department of Mathematics\\
100 Math Tower\\
Ohio State University\\
Columbus, OH 43210 (USA)}
\email{schoutens@math.ohio-state.edu}}

\newcommand{\emptyprop}{q}

\newcommand \id{\mathfrak a}

\newcommand \inverse[2]{{#1^{-1}(#2)}}
\newcommand \iso{\cong}

\newcommand \map[1]{{\newcommand{\tmpprop}{#1q}  \if\tmpprop\emptyprop \to\else \xrightarrow{{\phantom{i}{#1}\phantom{i}}}\fi}} 
\newcommand \maxim{\mathfrak m}
\newcommand \nat{\mathbb N}

\newcommand \pol[2]{#1[#2]}

\newcommand \pr{\mathfrak p}

\newcommand \rij[2]{(#1_1,\dots,#1_{#2})}

\newcommand \tensor{\otimes}
\newcommand \tor[4]{\operatorname{Tor}^{#1}_{#2}(#3,#4)}
\newcommand \op\operatorname

\newcommand{\commdiagram}[9][]{%
\begin{equation}
{\newcommand{\tmpprop}{#1q} 
\if\tmpprop\emptyprop \relax\else \label{#1}\fi}
\begin{aligned}%
\mbox{
\begin{picture}(130,90)%
\put(120,70){\vector( 0,-1){50}}%
\put(10,80){\vector( 1, 0){100}}%
\put(0,70){\vector( 0,-1){50}}%
\put(10,10){\vector( 1, 0){100}}%
\put(115,80){\makebox(0,0)[l]{$#4$}}%
\put(5,80){\makebox(0,0)[r]{$#2$}}%
\put(115,10){\makebox(0,0)[l]{$#9$}}%
\put(5,10){\makebox(0,0)[r]{$#7$}}%
\put(-3,50){\makebox(0,0)[r]{$#5$}}
\put(123,50){\makebox(0,0)[l]{$#6$}}
\put(60,3){\makebox(0,0)[c]{$#8$}}
\put(60,88){\makebox(0,0)[c]{$#3$}}
\end{picture}}
\end{aligned}
\end{equation}}

\newcommand\commtrianglefront[7][]{%
\begin{equation}
{\newcommand{\tmpprop}{#1q} 
\if\tmpprop\emptyprop \relax\else \label{#1}\fi}
\begin{aligned}%
\mbox{
\begin{picture}(120,80)%
\put(55,70){\vector(-1,-2){30}}
\put(65,70){\vector(1,-2){30}}
\put(30,5){\vector(1,0){60}}
\put(60,75){\makebox(0,0)[c]{$#2$}}
\put(25,5){\makebox(0,0)[r]{$#4$}}
\put(95,5){\makebox(0,0)[l]{$#6$}}
\put(60,0){\makebox(0,0)[c]{$#5$}}
\put(37,43){\makebox(0,0)[r]{$#3$}}
\put(83,43){\makebox(0,0)[l]{$#7$}}
\end{picture}}
\end{aligned}
\end{equation}}

\newcommand\commtriangleback[7][]{%
\begin{equation}
{\newcommand{\tmpprop}{#1q}
\if\tmpprop\emptyprop \relax\else \label{#1}\fi}
\begin{aligned}%
\mbox{
\begin{picture}(120,80)%
\put(55,70){\vector(-1,-2){30}}
\put(65,70){\vector(1,-2){30}}
\put(30,5){\vector(1,0){60}}
\put(60,75){\makebox(0,0)[c]{$#2$}}
\put(25,5){\makebox(0,0)[r]{$#6$}}
\put(95,5){\makebox(0,0)[l]{$#4$}}
\put(60,0){\makebox(0,0)[c]{$#5$}}
\put(37,43){\makebox(0,0)[r]{$#7$}}
\put(83,43){\makebox(0,0)[l]{$#3$}}
\end{picture}}
\end{aligned}
\end{equation}}

\newcommand{\name}[1]{{\sc#1}}
\newcommand \acf{algebraically closed field}

\newcommand \ch{characteristic}
\newcommand \homo{homomorphism}
\newcommand \CM{Coh\-en-Mac\-au\-lay}
\renewcommand\iff{if, and only if,}

\DeclareMathOperator*{\UP}{ulim}
\newcommand \up[1]{\UP_{#1\to\infty}}
\newcommand \ul[1]{\seq{#1}\infty}
\newcommand \seq[2]{#1\mathstrut_{#2}}
\newcommand \sr{approximation}
\newcommand \BS{Brian\c{c}on-Skoda}

\newcommand  \nsrat{non-standard difference rational}
\newcommand  \genrat{generically F-rational}
\newcommand  \genratity{generic F-rationality}
\newcommand  \genreg{weakly generically F-regular}
\newcommand  \genregglob{generically F-regular}
\newcommand  \genregity{weak generic F-regularity}
\newcommand  \los{\L os' Theorem}
\newcommand  \BCM[1]{\mathcal B(#1)}

\newcommand \frob[1]{\mathbf{F}_{#1}}
\newcommand \ulfrob{\frob\infty}

\newcommand \alg[1]{ {#1} ^{\text{alg}}}
\newcommand \ac[2]{{\mathbb{#1}_{#2}^{\text{alg}}}}
\newcommand \zet{\mathbb Z}

\title [Canonical Big Cohen-Macaulay Algebras]{Canonical Big Cohen-Macaulay Algebras and Rational Singularities}
\author{Hans Schoutens}
\thanks{Partially supported by a  grant from the National Science Foundation and by visiting positions at Paris VII and at the Ecole Normale Sup\'erieure.}
\date{28.02.2003}
\OSU

\begin{document}

\begin{abstract}   
We give a canonical construction of a balanced big \CM\ algebra for a domain of finite type over $\mathbb C$  by taking ultraproducts of absolute integral closures in positive \ch. This yields a new tight closure characterization of rational singularities in \ch\ zero.
\end{abstract}

\maketitle

\section{Introduction}

In \cite{HoHT}, \name{Hochster} proves the existence of big \CM\ modules for a large class of Noetherian rings containing a field. Recall that a module $M$ over a Noetherian  local ring $R$ is called a \emph{big \CM\ module}, if there is a system of parameters of $R$ which is $M$-regular (the adjective \emph{big} is used to emphasize that $M$ need not be finitely generated). He also exhibits in that paper the utility of big \CM\ modules in answering various homological questions. Often, one can even obtain a big \CM\ module $M$ such that \emph{every} system of parameters is $M$-regular; these are called \emph{balanced} big \CM\ modules. In \cite{HHbigCM}, \name{Hochster} and \name{Huneke} show that for equi\ch\ excellent local domains, one can even find a balanced big \CM\ \emph{algebra}, that is to say, $M$ admits the structure of a (commutative) $R$-algebra. In fact, for $R$ a local domain of positive \ch, they show that the absolute integral closure of $R$, denoted $R^+$, is  a (balanced) big \CM\ algebra (it is easy to see that this is false in \ch\ zero).  In \cite{HHbigCM2}, using lifting techniques  similar to the ones developed in the original paper of \name{Hochster}, they obtain also the existence of big \CM\ algebras in \ch\ zero. However, the construction is no longer canonical and one looses the additional information one had in positive \ch. Nonetheless, many useful applications follow, see \cite[\S 9]{HuTC} or \cite{HHbigCM2}.

In this paper, I will show  that for a  local domain $R$ of finite type over $\mathbb C$ (henceforth, a \emph{local $\mathbb C$-affine domain}), a simple construction of a balanced big \CM\ algebra  $\BCM R$ can be made, which restores canonicity, is weakly functorial and preserves many of the  good properties of the absolute integral closure. Namely, to the domain $R$, one associates certain \ch\ $p$ domains $\seq Rp$, called \emph{\sr{s}} of $R$, and of these one takes the absolute integral closure $\seq Rp^+$ and then forms the ultraproduct $\BCM R:=\up p\seq Rp^+$.  For generalities on ultraproducts, including \los, see \cite[\S 2]{SchNSTC}. Recall that an ultraproduct  of rings $\seq Cp$ is a certain homomorphic image of the direct product of the $\seq Cp$. This ultraproduct will be denoted by $\up p\seq Cp$, or simply by $\ul C$, and similarly, the image of a sequence $(\seq ap\mid p)$ in $\ul C$ will be denoted by $\up p\seq ap$, or simply by $\ul a$.

The notion of \sr\ goes back to the paper \cite{SchNSTC}, where it was introduced to define a closure operation, called \emph{non-standard tight closure} on $\mathbb C$-affine algebras by means of a so-called \emph{non-standard Frobenius}. Let me briefly recall the construction of an \sr\ (details and proofs can be found in \cite[\S 3]{SchNSTC}). Suppose $R$ is of the form $\pol{\mathbb C}X/I$, or possibly, a localization of such an algebra with respect to a prime ideal $\pr$. There is a fundamental (but non-canonical) isomorphism between the field of complex numbers on the one hand, and the ultraproduct of all the fields $\ac Fp$ on the other hand, where $\ac Fp$ denotes the algebraic closure of the $p$-element field. Therefore, for every element $c$ in $\mathbb C$, we can choose a representative in the product, that is to say,  a sequence of elements $\seq cp\in\ac Fp$, called an \emph\sr\ of $c$, such that $\up p\seq cp=c$. Applying this to each coefficient of a polynomial $f\in\pol{\mathbb C}X$ separately, we get a sequence of polynomials $\seq fp\in\pol{\ac Fp}X$ (of the same degree as $f$), called again an \emph{\sr} of $f$. If we apply this to the generators of $I$ and $\pr$, we generate ideals $\seq Ip$ and $\seq\pr p$ in $\pol{\ac Fp}X$, called once more \sr{s} of $I$ and $\pr$ respectively. One shows that $\seq\pr p$ is prime for almost all $p$. Finally, we set $\seq Rp:=\pol{\ac Fp}X/\seq Ip$ (or its localization at the prime ideal $\seq \pr p$) and call the collection of these \ch\ $p$ rings an \emph\sr\ of $R$. Although the choice of an \sr\ is not unique, almost all its members are the same; this is true for every type of \sr\ just introduced (here and elsewhere, \emph{almost all} means with respect to a non-specified but fixed non-principal ultrafilter). Moreover, if we depart from a different presentation of $R$ as a $\mathbb C$-affine algebra, then the resulting \sr\ is isomorphic to $\seq Rp$, for almost all $p$. In particular, the ultraproduct $\ul R:=\up p\seq Rp$ of the $\seq Rp$ is uniquely determined up to $R$-algebra isomorphism and is called the \emph{non-standard hull} of $R$. There is a natural embedding $R\to \ul R$, the main property of which was discovered by \name{van den Dries} in \cite{vdD79}: $R\to \ul R$  is faithfully flat (note that in general, $\ul R$ is no longer Noetherian nor even separated).  In case $R$ is a local domain, almost all $\seq Rp$ are local domains. Therefore, the ultraproduct $\BCM R$ of the $\seq Rp^+$ is well defined and unique up to $R$-algebra isomorphism and we get our first main result.

\begin{stelling}\label{T:A}
If $R$ is a local $\mathbb C$-affine domain, then $\BCM R$ is a balanced big \CM\ algebra.
\end{stelling}

In fact, due to canonicity, the operation of taking $\BCM\cdot$ is  weakly functorial (see Theorem~\ref{T:BCMeq} for a precise statement). Moreover, $\BCM R$ has the additional property that every monic polynomial over it splits completely in linear factors, so that $\BCM R$ is in particular Henselian. In $\BCM R$, any sum of prime ideals is either the unit ideal or else again a prime ideal. This is explained in Section~\ref{s:prop}. In Section~\ref{s:ratsing}, we relate the construction of $\BCM\cdot$ with generic tight closure (this is one of the alternative closure operations in \ch\ zero introduced in \cite{SchNSTC}). One immediate corollary of the canonicity of our construction is the following \ch\ zero version of the generalized \BS\ Theorem  in \cite[Theorem 7.1]{HHbigCM2}.

\begin{stelling}\label{T:BS}
If $R$ is a local $\mathbb C$-affine domain and $I$ an ideal of $R$ generated by $n$ elements, then the integral closure of $I^{n+k}$ is contained in $I^{k+1}\BCM R\cap R$, for every $k\in\nat$.
\end{stelling}

In \cite{SchNSTC} the same result is proven if we replace $I\BCM R\cap R$ by the generic tight closure of $I$. This suggests that the appropriate \ch\ zero equivalent of the conjecture that tight closure equals plus closure is the conjecture that $I\BCM R\cap R$ always equals the generic tight closure of $I$. We show that in any case, the former is contained in the latter. Moreover, we have equality for parameter ideals, that is to say, the \ch\ zero equivalent of \name{Smith}'s result in \cite{SmParId} also holds (for a further discussion, see \ref{s:Bcl}).  Using this, we give a characterization of rational singularities, in terms of generic tight closure, extending the results of \name{Hara} \cite{HaRat} and \name{Smith} \cite{SmFrat}, at least in the affine case. 

\begin{stelling}\label{T:RatSing}
If $R$ is a local $\mathbb C$-affine domain, then $R$ has rational singularities \iff\ there exists a system of parameters $\mathbf x$ such that  $\mathbf x\BCM R\cap R=\mathbf xR$.
\end{stelling}

Note that we need \name{Hara}'s result for the proof (see Theorem~\ref{T:Frat} for more details), which itself relies on some deep vanishing theorems. In \cite{SchLogTerm}, we will give a similar characterization for log-terminal singularities. Using the above results, we recover the \BS\ Theorem of \name{Lipman-Teissier}. Another application is a new proof of \name{Boutot}'s main result in \cite{Bou}, at least for Gorenstein rational singularities (this also generalizes the main result of \cite{SchRatSing}; for a further generalization, see \cite[Theorem B]{SchLogTerm}).

\begin{citedtheorem}[\name{Boutot} \cite{Bou}]
Let $R\to S$ be a (cyclically) pure \homo\ of local $\mathbb C$-affine algebras. If $S$ is Gorenstein and has rational singularities, then $R$ has rational singularities. 
\end{citedtheorem}

In  the final section, some results of \cite{SchBetti} are extended to the present \ch\ zero situation. In particular, we obtain the following regularity criterion (see Theorem~\ref{T:reg}).

\begin{stelling}\label{T:isosing}
Let $R$ be a local $\mathbb C$-affine domain $R$ with  residue field $k$. If $R$ has an isolated singularity and  $\tor R1{\BCM R}k=0$, then $R$ is regular.
\end{stelling}

In contrast with the prime \ch\ case, I do not know whether for arbitrary local $\mathbb C$-affine domains $R$, the flatness of $R\to\BCM R$ is equivalent with the regularity of $R$.

\subsection*{Remark on the base field}
To make the exposition more transparent, I have only dealt in the text with the case that the base field is $\mathbb C$. However, the results extend to arbitrary uncountable base fields of \ch\ zero by the following observations. First, any uncountable \acf\ of \ch\ zero is the ultraproduct of (algebraically closed) fields of positive \ch\ by the Lefschetz Principle (see for instance \cite[Remark 2.5]{SchNSTC}) and this is the only property we used of $\mathbb C$. Second, if $A$ is a local $K$-affine domain with $K$ an arbitrary uncountable field, then $A^+$ is a $K^{\text{alg}}$-algebra, where  $K^{\text{alg}}$ is the algebraic closure of $K$. Therefore, in order to define $\BCM A$ in case $K$ has moreover \ch\ zero, we may replace $A$ by $A\tensor_K  K^{\text{alg}}$ and assume form the start that $K$ is uncountable and algebraically closed, so that   our first observation applies.

In a future paper, I will discuss how one can extend the quasi-hull $\BCM \cdot$ to equi\ch\ complete local domains. In \cite{SchMixBCMCR,SchMixBCM} we use the same techniques to obtain an asymptotic version of big \CM\ algebras in mixed \ch.

\section{Big \CM\ algebras}

\subsection{Absolute Integral Closure}
Let $A$ be a domain. The \emph{absolute integral closure} $A^+$ of $A$ is defined as follows. Let $Q$ be the  field of fractions of $A$ and let $\alg Q$ be its algebraic closure, We let $A^+$ be the integral closure of $A$ in $\alg Q$. Since algebraic closure is unique up to isomorphism, any two absolute integral closures of $A$ are isomorphic as $A$-algebras. To not have to deal with exceptional cases separately, we put $A^+=0$ if $A$ is not a domain.

In this paper, we will use the term \emph{$K$-affine algebra} for an algebra of finite type over a field $K$ or a localization of such an algebra with respect to a prime ideal; the latter will also be referred to as a \emph{local $K$-affine algebra}.

\subsection{Approximations and non-standard hulls}\label{s:nsh}
Let $A$ be a $\mathbb C$-affine algebra and choose an \sr\ $\seq Ap$ of $A$ (see the introduction;  for a precise definition and proofs, see \cite[\S 3]{SchNSTC}). The ultraproduct of the $\seq Ap$ is called the \emph{non-standard hull} of $A$ and is often denoted $\ul A$. The assignment $A\mapsto \ul A$ is functorial.  There is a natural \homo\ $A\to \ul A$, which is faithfully flat by \cite[Theorem 1.7]{SvdD}  (for an alternative proof, see \cite[A.2]{SchBArt}). It follows that if $I$ is an ideal in $A$ and $\seq Ip$ an \sr\ of $I$, then $I\ul A$ is the ultraproduct of the $\seq Ip$ and $I=I\ul A\cap A$. By \cite[Theorem 4.4]{SchNSTC}, almost all $\seq Ap$ are domains (respectively, local) \iff\ $A$ is a domain (respectively, local) \iff\ $\ul A$ is a domain (respectively, local). 

\subsection{The quasi-hull $\BCM\cdot$}\label{s:BCM}
Let $A$ be a $\mathbb C$-affine domain with \sr\ $\seq Ap$.  Define $\BCM A$ as the ultraproduct
	\begin{equation*}
	\BCM A:=\up p \seq Ap^+.
	\end{equation*}
In view of the uniqueness of the absolute integral closure,  $\BCM A$ is independent of the choice of the $\seq Ap$ and hence is uniquely determined by $A$ up to $A$-algebra isomorphism. If $A$ is local, then so is $\BCM A$. Given a \homo\ $A\to B$ of $\mathbb C$-affine algebras, we obtain \homo{s} $\seq Ap\to \seq Bp$, for almost all $p$, where $\seq Bp$ is an \sr\ of $B$ (see \cite[3.2.4]{SchNSTC}).  These \homo{s} induce (non-canonically) \homo{s} $\seq Ap^+\to\seq Bp^+$, which, in the ultraproduct,  yield a \homo\ $\BCM A\to\BCM B$.

Note that the natural \homo\ $A\to\BCM A$ factors through the non-standard hull $\ul A$, and in particular, $A\to \BCM A$ is no longer integral. Using \los\ and results on the absolute integral closure in \cite{HHbigCM} (see also \cite[Chapter 9] {HuTC}), we get the following more precise version of Theorem~\ref{T:A}.

\begin{theorem}\label{T:BCMeq}
For each local $\mathbb C$-affine domain $A$, the $A$-algebra $\BCM A$ is a \emph{balanced big \CM\ algebra} in the sense that any system of parameters $\mathbf x$ of $A$ is a $\BCM A$-regular sequence. Moreover, if $A\to B$ is a \homo\  between local  $\mathbb C$-affine domains, then there exists a $\mathbb C$-algebra \homo\ $\BCM A\to\BCM B$ giving rise to a commutative diagram
	\commdiagram [bcm] A {} B {} {} {\BCM A} {} {\BCM B.}

If $A\to B$ is  finite and injective, then $\BCM A=\BCM B$.
\end{theorem}
\begin{proof}
Let $\seq Ap$  be an \sr\ of $A$ and $\ul A$ its non-standard hull. Let $\mathbf x$ be a system of parameters in $A$ with \sr\ $\seq {\mathbf x}p$.  By \cite[Theorem 4.5]{SchNSTC} almost all $\seq{\mathbf x}p$ are a system of parameters of $\seq Ap$. Therefore, by  \cite[Theorem 1.1]{HHbigCM}, the sequence $\seq{\mathbf x}p$ is $\seq Ap^+$-regular, for almost all $p$. \los\ then yields that $\mathbf x$ is a $\BCM A$-regular sequence.

The existence of the \homo\ $\BCM A\to \BCM B$ and the commutativity of diagram~\eqref{bcm} follow from the above discussion. Finally, if $B$ is finite overring of $A$, then by \cite[Theorem 4.7]{SchNSTC}, so will almost all $\seq Bp$ be over $\seq Ap$, where $\seq Ap$ and $\seq Bp$ are \sr{s} of $A$ and $B$ respectively. In particular, $\seq Ap^+=\seq Bp^+$, for almost all $p$, proving that $\BCM A=\BCM B$.
\end{proof}

\begin{corollary}\label{C:reg}
For each local $\mathbb C$-affine regular ring $A$, the natural map $A\to \BCM A$ is faithfully flat.
\end{corollary}
\begin{proof}
It is well-known that a balanced big \CM\ module over a regular local ring is flat (see for instance \cite[Theorem IV.1]{SchFPD} or \cite[Lemma 2.1(d)]{HHbigCM2}). Since the natural map $A\to\BCM A$ is local, the result follows.
\end{proof}

As in positive \ch, we can construct big \CM\ algebras over any reduced local $\mathbb C$-affine ring $A$ by letting $\BCM A$ be the product of all $\BCM{A/\pr}$, where $\pr$ runs over all minimal prime ideals of $A$. As for localization, we have a slightly less pretty result as in positive \ch: if $A$ is a local $\mathbb C$-affine domain with non-standard hull $\ul A$  and if $\pr$ is a prime ideal of $A$, then 
	\begin{equation}\label{eq:loc}
	\BCM {A_\pr}\iso \BCM A\tensor_{\ul A}(\ul A)_{\pr\ul A}.
	\end{equation}
Indeed, if $\seq Ap$ and $\seq\pr p$ are \sr{s} of $A$ and $\pr$ respectively, then by \cite[Lemma 6.5]{HHbigCM}, we have isomorphisms 
	\begin{equation*}
	((\seq Ap)_{\seq\pr p})^+\iso (\seq Ap^+)_{\seq\pr p} = \seq Ap^+\tensor_{\seq Ap}(\seq Ap)_{\seq\pr p}. 
	\end{equation*}
Taking ultraproducts, we get isomorphism \eqref{eq:loc}. It follows that Corollary~\ref{C:reg} also holds if we drop the requirement that $A$ is local (use that $\BCM A_\pr\to\BCM {A_\pr}$ is faithfully flat, for every prime ideal $\pr$ of $A$, by \eqref{eq:loc}). We also obtain the following \ch\ zero analogue of \cite[Theorem 6.6]{HHbigCM}.

\begin{theorem}
If $A$ is a $\mathbb C$-affine domain and  $I$ an ideal in $A$ of height $h$, then $H^j_I(\BCM A)=0$, for all $j<h$.
\end{theorem}
\begin{proof}
As in the proof of \cite[Theorem 6.6]{HHbigCM}, it suffices to show that for every maximal ideal $\maxim$ of $A$ containing $I$, we have that $H^j_I(\BCM A)_\maxim=0$, for $j<h$. Since $A_\maxim\to (\ul A)_{\maxim\ul A}$ is faithfully flat, as explained in Section~\ref{s:nsh}, it suffices to show that 
	\begin{equation*}
	H^j_I(\BCM A)\tensor_{\ul A}(\ul A)_{\maxim\ul A}=0.
	\end{equation*}
 By  \eqref{eq:loc}, the left hand side is simply $H^j_I(\BCM{A_\maxim})$ and therefore, the problem reduces to the case that $A$ is local. Let   $\rij xh$ be part of a system of parameters of $A$  contained in $I$. Since $\rij xh$ is $\BCM A$-regular by Theorem~\ref{T:BCMeq}, the vanishing of $H^j_I(\BCM A)$ for $j<h$ is then clear since local cohomology can be viewed as a direct limit of Koszul cohomology.
\end{proof}

\section{Properties of $\BCM A$}\label{s:prop}

Let us call a domain $S$ \emph{absolutely integrally closed} if every monic polynomial over $S$ has a root in $S$.

\begin{lemma}\label{L:aic}
 For a domain $S$ with field of fractions $Q$, the following are equivalent.
\begin{enumerate}
\item\label{i:aic} $S$ is absolutely integrally closed.
\item\label{i:mon} Every monic polynomial completely splits in $S$.
\item\label{i:ac} $S$ is integrally closed in $Q$ and $Q$ is algebraically closed.
\end{enumerate}
\end{lemma}
\begin{proof}
The implications $\eqref{i:ac}\implies \eqref{i:mon}$ and $\eqref{i:mon}\implies \eqref{i:aic}$ are straightforward. Hence assume that $S$ is absolutely integrally closed. It is clear that $S$ is then integrally closed in $Q$. So remains to show that $Q$ is algebraically closed. In other words, we have to show that every non-zero one-variable polynomial $F\in\pol QT$ has a root in $Q$. Clearing denominators, we may assume that $F\in\pol ST$. Let $a\in S$ be the (non-zero) leading coefficient of $F$ and $d$ its degree. We can find a monic polynomial $G$ over $S$, such that $a^{d-1}F(T)=G(aT)$. By assumption, $G(b)=0$ for some $b\in S$. Hence $F(b/a)=0$, as required.
\end{proof}

It follows from \cite[Lemma 6.5]{HHbigCM} that a domain $S$ is the absolute integral closure of a subring $A$ \iff\ $S$ is  absolutely integrally closed and $A\subset S$ is integral.

\begin{proposition}\label{P:mon}
If $A$ is a $\mathbb C$-affine domain, then $\BCM A$ is absolutely integrally closed. In particular, if $A$ is local, then $\BCM A$ is Henselian.
\end{proposition}
\begin{proof}
Let $F(T):=T^d+a_1T^{d-1}+\dots+a_d$  be a monic polynomial in the single variable $T$ with $a_i\in\BCM A$. We need to show that $F$ has a root in $\BCM A$. Choose $\seq{a_i}p\in\seq Ap^+$, such that $\up p\seq{a_i}p=a_i$, for all $i$, where $\seq Ap$ is some \sr\ of $A$. Hence we can find $\seq{b}p\in\seq Ap^+$ such that
	\begin{equation*}
	(\seq bp)^d+\seq{a_1}p(\seq bp)^{d-1}+\dots+\seq{a_d}p= 0.
	\end{equation*}
Therefore, by \los, $b:=\up p\seq{b}p$ is a root of $F$. The last statement is then immediate by definition of Henselian.
\end{proof}

\begin{corollary}
Let $A$ be a $\mathbb C$-affine domain. The sum of any collection of prime ideals in $\BCM A$ is either prime or the unit ideal. If $\mathfrak g_i$ are $\pr_i$-primary ideals, for $i$ in some index set $I$, and if $\pr:=\sum_{i\in I}\pr_i$ is not the unit ideal, then $\sum_{i\in I}\mathfrak g_i$ is $\pr$-primary.
\end{corollary}
\begin{proof}
By Proposition~\ref{P:mon}, the ring $\BCM A$ is quadratically closed and therefore has the stated properties by \cite[Theorem 9.2]{HHbigCM}.
\end{proof}

The next result shows that $\BCM A$, viewed as an $\ul A$-module, also behaves very much like a \CM\ module.

\begin{proposition}\label{P:CM}
Let $A$ be a local $\mathbb C$-affine domain. Let $\rij xd$ be part of a system of parameters of $A$  and let $\pr_1,\dots,\pr_s$ be the minimal prime ideals of $\rij xdA$. If $\ul t\in \ul A$ does not lie in any $\pr_i\ul A$, then $(x_1,\dots,x_d,\ul t)$ is a $\BCM A$-regular sequence. 
\end{proposition}
\begin{proof}
Suppose $\ul t\in \ul A$ lies outside all $\pr_i\ul A$ and suppose $\ul b\in \BCM A$ is such that 
	\begin{equation*}
	\ul t\ul b\in \rij xd\BCM A.
	\end{equation*} 
Let $\seq Ap$, $\seq {x_i}p$ and $\seq{\pr_i} p$ be \sr{s} of $A$, $x_i$ and $\pr_i$ respectively. It follows from \cite[Theorem 4.5]{SchNSTC} that $(\seq{x_1}p,\dots,\seq{x_d}p)$ is part of a system of parameters in $\seq Ap$, and from \cite[Theoprem 4.4]{SchNSTC}, that $\seq{\pr_1} p,\dots,\seq{\pr_s}p$ are the minimal prime ideals of $(\seq{x_1}p,\dots,\seq{x_d}p)\seq Ap$, for almost all $p$. Choose $\seq tp$ and $\seq bp$ in $\seq Ap$ and $\seq Ap^+$ respectively such that their ultraproduct is $\ul t$ and $\ul b$. By \los, almost all $\seq tp$ lie outside any $\seq{\pr_i} p$, and $\seq tp\seq bp\in(\seq{x_1}p,\dots,\seq{x_d}p)\seq Ap^+$. Therefore, $(\seq{x_1}p,\dots,\seq{x_d}p,\seq tp)$ is part of a system of parameters in $\seq Ap$ and hence, by  \cite[Theorem 1.1]{HHbigCM},  is an $\seq Ap^+$-regular sequence, for almost all $p$. It follows that $\seq bp\in(\seq{x_1}p,\dots,\seq{x_d}p)\seq Ap^+$, for almost all $p$, whence, by \los, that $\ul b\in \rij xd\BCM A$.
\end{proof}

\section{Rational Singularities}\label{s:ratsing}

Recall the definition of generic tight closure from  \cite{SchNSTC}. Let $A$ be a (local) $\mathbb C$-affine  algebra, $I$ an ideal of $A$ and $z$ an arbitrary element. We say that $z$ lies in the \emph{generic tight closure} of   $I$, if $\seq zp$ lies in the tight closure of $\seq Ip$, for almost all $p$, where $\seq zp$ and $\seq Ip$ are some \sr{s} of $z$ and $I$ respectively. In \cite{SchNSTC} it is shown that this yields a closure operation with similar properties as \ch\ zero tight closure, and that it is contained in non-standard tight closure (for the definition of non-standard (tight) closure  and for further properties of these closure operations,  see \cite{SchNSTC}; variants can be found in \cite{SchHR,SchRatSing}). 

\begin{corollary}\label{C:parid}
Let $R$ be a local $\mathbb C$-affine domain and let $I$ be an ideal generated by a system of parameters of $R$. The generic tight closure of $I$ is equal to $I\BCM R\cap R$.

More generally, for arbitrary $I$, we have that $I\BCM R\cap R$ is contained in the generic tight closure of $I$ (whence in the non-standard tight closure of $I$).
\end{corollary}
\begin{proof}
Let $\seq Rp$ and $\seq Ip$ be \sr{s} of $R$ and $I$ respectively. Let $f\in R$ with \sr\ $\seq fp$. Assume first that $f\in I\BCM R$. It follows that $\seq fp\in\seq Ip\seq Rp^+$, for almost all $p$. Since in general, $JB\cap A$ lies in the tight closure of $J$, for any integral extension $A\to B$ of prime \ch\ rings and any ideal $J\subset A$ (\cite[Theorem 1.7]{HuTC}), we get that $\seq fp$ lies in the tight closure of $\seq Ip$, for almost all $p$. However, this just means that $f$ lies in the generic tight closure of $I$. Conversely, if $f$ lies in the generic tight closure of $\mathbf xR$, where $\mathbf x$ is a system of parameters with \sr\ $\seq{\mathbf x}p$, then $\seq fp$ lies in the tight closure of $\seq {\mathbf x}p\seq Rp$ and $\seq {\mathbf x}p$ is a system of parameters in $\seq Rp$  by \cite[Theorem 4.5]{SchNSTC}, for almost all $p$. By the  result of \name{Smith} in \cite{SmParId}, tight closure equals `plus closure' for any ideal generated by a system of parameters, so that $\seq fp\in\seq {\mathbf x}p\seq Rp^+$. Taking ultraproducts, we get that $f\in \mathbf x\BCM R$.
\end{proof}

From this it is clear that Theorem~\ref{T:BS} is a strengthening of the \BS\ Theorem in \cite{SchNSTC} (see also \cite{SchBS}).  We also get the following sharpening of \cite[Theorem 6.2]{SchRatSing} (its converse also holds and will be proved in Theorem~\ref{T:Frat} below). 

\begin{theorem}\label{T:ratsing}
If   a local $\mathbb C$-affine domain $R$ admits a  system of parameters  $\mathbf x$  such that  $\mathbf xR=\mathbf x\BCM R\cap R$, then $R$ has rational singularities.
\end{theorem}
\begin{proof}
Let $\mathbf x:=\rij xd$. Let us first show that 
	\begin{equation}\label{eq:parid}
	\rij xiR=\rij xi\BCM R\cap R,
	\end{equation}
for all $i$. Let $I_i$ denote the ideal $\rij xiR$ and put $J_i:=I_i\BCM R\cap R$. We want to show that $I_i=J_i$, for all $i$, and we will achieve this by a downward induction on $i$.  The case $i=d$ holds by assumption. Suppose we showed already that $I_{i+1}=J_{i+1}$. Let $z\in J_i$. In particular, $z\in J_{i+1}=I_{i+1}$, so that we can write $z=a+rx_{i+1}$, for some $a\in I_i$ and some $r\in R$. Hence $z-a=rx_{i+1}\in J_i\subset I_i\BCM R$. Since $x_{i+1}$ is a non-zero divisor modulo $I_i\BCM R$ by Theorem~\ref{T:BCMeq}, we get that $r\in I_i\BCM R$, whence $r\in J_i$. In conclusion, we showed that $J_i=I_i+x_{i+1}J_i$. Nakayama's Lemma therefore yields $I_i=J_i$. 

Next, we show that $\mathbf x$ is $R$-regular. Suppose $zx_{i+1}\in I_i$. Since $\rij xd$ is $\BCM R$-regular by  Theorem~\ref{T:BCMeq}, we get $z\in I_i\BCM R$. By \eqref{eq:parid}, we therefore have $z\in I_i$, showing that $\mathbf x$ is $R$-regular. It follows that $R$ is \CM. By Corollary~\ref{C:parid} and  \eqref{eq:parid} also every principal height one ideal is equal to its generic tight closure. By \cite[Theorem 4.6 and Remark 4.7]{SchRatSing}, this implies that $R$ is normal. Finally, using the fact that $I_d$ is equal to its own generic tight closure by Corollary~\ref{C:parid}, we can repeat the argument in the proof of \cite[Theorem 6.2]{SchRatSing} to conclude that $R$ has rational singularities (see \cite[Remark 6.3]{SchRatSing}).
\end{proof}

\begin{definition}
Call a local $\mathbb C$-affine domain \emph{\genrat}, if  some ideal generated by a system of parameters is equal to its own generic tight closure. If, in contrast,  every ideal is equal to its own generic tight closure, then we will call such a domain \emph\genreg.

Similarly, for $R$ a local $\mathbb C$-affine domain, we say that $R$ is \emph{$\mathcal B$-rational}, if $\mathbf x\BCM R\cap R=\mathbf xR$, for some system of parameters $\mathbf x$ of $R$. If $R\to \BCM R$ is cyclically pure (that is to say, $I=I\BCM R\cap R$, for every ideal $I$ of $R$), we say that $R$  is \emph{weakly $\mathcal B$-regular}. 

If every localization of $R$ at a prime ideal is \genreg\ (respectively, weakly $\mathcal B$-regular), then we call $R$ \emph{\genregglob} (respectively, $\mathcal B$-regular).
\end{definition}

With this terminology, Theorem~\ref{T:ratsing} shows that a $\mathcal B$-rational local $\mathbb C$-affine domain has rational singularities (see also Theorem~\ref{T:Frat} and Definition~\ref{s:Bcl} below). The notion of (weak) $\mathcal B$-regularity is reminiscent of the notion \emph{(weak) CM$\mathstrut^n$-regularity} from \cite{HHbigCM2}. Conjecturally, \genreg\ and \genregglob\ are equivalent, and so are their $\mathcal B$-analogues expected to be (in positive \ch\ the latter holds automatically, but not so in the present case, due to the more complicated nature of the localization of $\BCM R$ given by \eqref{eq:loc}).  Surprisingly,  Corollary~\ref{C:parid} not only yields that  \genreg\ implies weakly $\mathcal B$-regular, but even $\mathcal B$-regular.

\begin{corollary}\label{C:Breg}
If a local $\mathbb C$-affine domain is \genreg, then it is $\mathcal B$-regular. 
\end{corollary}
\begin{proof}
Let $R$ be a \genreg\ local $\mathbb C$-affine domain, $I$ an arbitrary ideal in $R$ and $\pr$ a prime ideal. Put $S:=R_\pr$. We need to show that any $a$ in $I\BCM{S}\cap S$ lies already in $IS$. Let $\seq Rp$, $\seq ap$, $\seq Ip$ and $\seq\pr p$ be \sr{s} of $R$, $a$, $I$ and $\pr$ respectively, and put $\seq Sp:=(\seq Rp)_{\seq\pr p}$. By \los,  
	\begin{equation*}
	\seq ap \in \seq Ip \seq Sp^+=\seq Ip (\seq Rp^+)_{\seq\pr p}
	\end{equation*}
for almost all $p$, where the equality follows from \cite[Lemma 6.5]{HHbigCM}. Hence there exists $\seq sp\in\seq Rp$ but not in $\seq\pr p$ so that $\seq sp\seq ap$ lies in $\seq Ip\seq Rp^+$. Since plus closure is contained in tight closure, we get that $\seq sp\seq ap$ lies in the tight closure of $\seq Ip$. By assumption, the latter is tightly closed for almost all $p$, so that $\seq sp\seq ap\in\seq Ip$ whence $\seq ap\in \seq Ip\seq Sp$. Taking ultraproducts, we get that $a\in I\ul S$. Since $S\to \ul S$ is faithfully flat, we get that $a\in IS$, as required.
\end{proof}

\begin{proposition}\label{P:genrat}
For $R$ a local $\mathbb C$-affine domain with \sr\ $\seq Rp$, almost all $\seq Rp$ are F-rational \iff\ $R$ is \genrat.
\end{proposition}
\begin{proof}
Let $\mathbf x$ be a system of parameters of $R$ and let $\seq {\mathbf x}p$ be an \sr\ of $\mathbf x$. By \cite[Theorem 4.5]{SchNSTC} almost all $\seq {\mathbf x}p$ are a system of parameters of $\seq Rp$. Suppose first that almost all  $\seq Rp$ are F-rational. Let $y$ be in the generic tight closure of $\mathbf xR$ and let $\seq yp$ be an \sr\ of $y$. Hence almost all $\seq yp$ lie  in the tight closure of $\seq {\mathbf x}p\seq Rp$, whence in $\seq {\mathbf x}p\seq Rp$  by F-rationality. Therefore, $y\in \mathbf xR$ by \los.

Conversely, assume almost all $\seq Rp$ are not F-rational. This means that for almost all $p$, the tight closure of $\seq {\mathbf x}p\seq Rp$ is strictly bigger than $\seq {\mathbf x}p\seq Rp$. Let $\ul J$ be the ultraproduct of the  tight closures of the $\seq {\mathbf x}p\seq Rp$ . By \los, $\mathbf x\ul R\varsubsetneq \ul J$. Since $\mathbf xR$ is primary to the maximal ideal in $R$, we have an isomorphism $R/\mathbf xR\iso \ul R/\mathbf x\ul R$ (use for instance \cite[Theorem 4.5]{SchNSTC}). Symbolically, this means that   $\ul R=R+\mathbf x\ul R$ (as sets), and hence that $\ul J= (\ul J\cap R)+\mathbf x\ul R$. Therefore, putting $J:=\ul J\cap R$, we showed that $\ul J=J\ul R$. Since $\mathbf x\ul R\varsubsetneq \ul J$, we get that $\mathbf xR\varsubsetneq J$. However, one easily checks that $J$ is just the generic tight closure of $\mathbf xR$. Hence, for no system of parameters $\mathbf x$, is $\mathbf xR$ equal to its generic tight closure, showing that $R$ is not \genrat.
\end{proof}

\begin{remark}\label{R:genreg}
In the course of the proof we actually established the following more general result. Let $(R,\maxim)$ be a local $\mathbb C$-affine domain and let $I$ be $\maxim$-primary. The ultraproduct of the tight closures of an \sr\ of $I$ is equal to the extension of the generic tight closure of $I$ to $\ul R$. It follows that,  if almost all $\seq Rp$ are weakly F-regular, then $R$ is \genreg. Indeed, let $\tilde I$ be the generic tight closure of an ideal $I$ and let $\seq Ip$ be an \sr\ of $I$. Suppose first that $I$ is $\maxim$-primary. By what we just said, $\tilde I\ul R$ is then equal to the ultraproduct of the $(\seq Ip)^* =\seq Ip$, that is to say, equal to $I\ul R$. Hence by faithful flatness, $I=\tilde I$. For $I$ arbitrary, $\tilde I$ is contained in the generic tight closure of $I+\maxim^n$, and by the previous argument that is just $I+\maxim^n$. Since this holds for all $n$,  Krull's Intersection Theorem yields $I=\tilde I$. 

However, this argument does not prove the converse (since the ideals that disprove the weak F-regularity of each $\seq Rp$ might be of unbounded degree). Nonetheless, we suspect the converse to be true as well. Proposition~\ref{P:Gor} below gives the converse under the additional Gorenstein assumption.
\end{remark}

\begin{proposition}\label{P:regseq}
For  a  local $\mathbb C$-affine domain $R$, the following are true.
\begin{enumerate}
\item\label{F:regseq}
If $\rij xd$ is a regular sequence for which 
	\begin{equation*}
	\rij xd\BCM R\cap R=\rij xdR, 
	\end{equation*}
then 
	\begin{equation*}
	(x_1^t,\dots,x_d^t)\BCM R\cap R= (x_1^t,\dots,x_d^t)R, 
	\end{equation*}
for all $t\geq 1$.
\item\label{F:colon}
If $I$ is an ideal of $R$ for which $I\BCM R\cap R=I$ and if $J$ is an arbitrary ideal  of $R$, then $(I:J)\BCM R\cap R=(I:J)$.
\item\label{F:par}
If $R$ is $\mathcal B$-rational, then $I\BCM R\cap R=I$, for every ideal $I$ generated by part of a system of parameters.
\end{enumerate}
\end{proposition}
\begin{proof}
We translate the usual tight closure proofs from \cite{HuTC} to the present situation. For \eqref{F:regseq},  induct on $t$, where $t=1$ is just the hypothesis. Let $z$ be an element in 
	\begin{equation*}
	(x_1^t,\dots,x_d^t)\BCM R\cap R.
	\end{equation*} 
If $x_iz\notin (x_1^t,\dots,x_d^t)R$, then we may replace $z$ by $x_iz$. Therefore, we may assume without loss of generality that $zI\subset (x_1^t,\dots,x_d^t)R$. Since $\rij xd$ is $R$-regular, $((x_1^t,\dots,x_d^t)R:I)=(x_1^t,\dots,x_d^t,x^{t-1})R$, where $x$ is the product of all $x_i$. Hence we may write $z=wx^{t-1}$, for some $w\in R$. By assumption, $z=wx^{t-1}\in (x_1^t,\dots,x_d^t)\BCM R$. Since $\rij xd$ is $\BCM R$-regular by Theorem~\ref{T:BCMeq}, we get that $w\in I\BCM R$, whence $w\in I$, by the original hypothesis. However, this shows that $z=wx^{t-1}\in (x_1^t,\dots,x_d^t)R$.

Assertion~\eqref{F:colon} is clear, since $z\in (I:J)\BCM R\cap R$ implies that $zJ\subset I\BCM R\cap R=I$. To prove the last assertion, assume that $R$ is $\mathcal B$-rational, say, $\mathbf x\BCM R\cap R=\mathbf x R$ for some  system of parameters $\mathbf x+\rij xd$. Let $I$ be an ideal generated by an arbitrary system of parameters $\rij yd$. Since we can calculate the top local cohomology group $\op H_\maxim^d(R)$ as the direct limit of the system $R/(x_1^t,\dots,x_d^t)R$ or, alternatively, as the direct limit of the system $R/(y_1^t,\dots,y_d^t)R$, we must have an embedding $R/I\to R/(x_1^t,\dots,x_d^t)R$ for sufficiently large $t$. Put differently, for large enough $t$, we have that $I=((x_1^t,\dots,x_d^t)R: a)$, for some $a\in R$ (see for instance \cite[Exercise 4.4]{HuTC}). It follows therefore from \eqref{F:regseq} and \eqref{F:colon} that $I\BCM R\cap R=I$. If $I$ is only generated by part of a system of parameters, then the assertion follows from \eqref{eq:parid} in the proof of Theorem~\ref{T:ratsing}.
\end{proof}

By virtually the same argument, assertion~\eqref{F:par} also holds for \genrat\ rings.

\subsection{Models}
Let $K$ be a field and $R$ a $K$-affine algebra. With a \emph{model of $R$} (called \emph{descent data} in \cite{HHZero}) we mean a pair $(Z,R_Z)$ consisting of a subring $Z$ of $K$ which is finitely generated over $\zet$ and a $Z$-algebra $R_Z$ essentially of finite type, such that $R\iso R_Z\tensor_Z K$. Oftentimes,  we will think of $R_Z$  as being the model. Clearly, the collection of models $R_Z$ of $R$ forms a direct system whose union is $R$. We say that $R$ has \emph{F-rational type} (respectively, has \emph{weakly F-regular type}), if there  exists a model $(Z,R_Z)$, such that $R_Z/\pr R_Z$ is F-rational (respectively, weakly F-regular) for all maximal ideals $\pr$ of $Z$ (note that we may always localize $Z$ at a suitably chosen element so that the property holds for all maximal ideals). See  \cite{HHZero} or \cite{HuTC} for more details.

In order to compare the notions of F-rational type and \genratity, we need to better understand the relation between reduction modulo $p$ and \sr{s}. We will see that \sr{s} are base changes to the algebraic closure of the residue field of reductions modulo $p$, where the choice of the embedding of the residue field in its algebraic closure is determined by the ultrafilter.

\begin{lemma}\label{L:srzet}
Let $Z$ be a finitely generated $\zet$-subalgebra of $\mathbb C$. For almost all $p$, there exists a \homo\ $\seq\gamma p\colon Z\to \ac Fp$, such that the sequence $\seq\gamma p(z)$ is an \sr\ of $z$, for each $z\in Z$.
\end{lemma}
\begin{proof}
Write $Z\iso\pol\zet Y/\rij gm\pol\zet Y$, with $Y$ a finite tuple of variables. Let $\mathbf y$ be the image of the tuple $Y$ in $\mathbb C$ under the embedding $Z\subset \mathbb C$ and take an \sr\ $\seq{\mathbf y}p$ of $\mathbf y$ in $\ac Fp$. By \los,  $\rij gm\pol {\mathbb F_p}Y$ is contained in the kernel of the algebra \homo\ $\pol {\mathbb F_p}Y\to \ac Fp$ given by $Y\mapsto \seq{\mathbf y}p$, for almost all $p$. This induces a \homo\ $\seq \gamma p\colon Z\to \ac Fp$ as asserted. Remains to verify the \sr\ property. To this end, let $z\in Z$ be represented by the image of $G\in\pol\zet Y$, that is to say, $z=G(\mathbf y)$. By construction, $\seq\gamma p(z)=G(\seq{\mathbf y}p)$. Since in the ultraproduct
	\begin{equation*}
	\up p G(\seq{\mathbf y}p) =  G(\up p\seq{\mathbf y}p) =G(\mathbf y) =z,
	\end{equation*}
we showed that $\seq\gamma p(z)$ is an \sr\ of $z$.
\end{proof}

Note that almost all $\seq\gamma p(Z)\subset \ac Fp$ are in fact separable field extensions.

\begin{corollary}\label{C:model}
Let $R$ be a local $\mathbb C$-affine domain with \sr\ $\seq Rp$. For each finite subset of $R$, we can find a   model $(Z,R_Z)$ of $R$ containing this subset, and, for almost all $p$,  a \homo\ $\seq\gamma p\colon Z\to \ac Fp$ inducing a separable field extension $\seq\gamma p(Z)\subset \ac Fp$, such that 
	\begin{equation}\label{eq:seqrp}
	\seq Rp:=R_Z \tensor_Z \ac Fp
	\end{equation}
is an \sr\ of $R$.  

Moreover, for each $r\in R_Z$, we get an \sr\ of $r$ by taking its image in $\seq Rp$ via the canonical \homo\ $R_Z\to \seq Rp$.
\end{corollary}
\begin{proof}
Suppose $R$ is the localization of $\pol{\mathbb C}X/I$ at the prime ideal $\maxim$. Take any model $(Z,R_Z)$ of $R$ containing the prescribed subset. After possibly enlarging this model, we may moreover assume that there exists ideals $I_Z$ and $\maxim_Z$ in $\pol ZX$ such that 
	\begin{equation*}
	R_Z= (\pol ZX/I_Z)_{\maxim_Z}
	\end{equation*}
 (whence $I=I_Z\pol{\mathbb C}X$ and $\maxim=\maxim_Z\pol{\mathbb C}X$).  Let $\seq\gamma p\colon Z\to \ac Fp$ be a \homo\ as in Lemma~\ref{L:srzet} such that  $\seq\gamma p(z)$ is an \sr\ of $z$, for each $z\in Z$. Let $\seq Ip$ (respectively,  $\seq \maxim p$) be the ideal  in $\pol{\ac Fp}X$ generated by all $f^{\seq\gamma p}$ with $f\in I_Z$ (respectively, $f\in\maxim_Z$), where we write $f^{\seq\gamma p}$  for the polynomial obtained from $f$ by applying $\seq\gamma p$ to each of its coefficients. It follows that $\seq Ip$ and $\seq \maxim p$) are \sr{s} of $I$ and $\maxim$ respectively. Therefore
	\begin{equation*}
	(\pol{\ac Fp} X/\seq Ip)_{\seq\maxim p}\iso R_Z\tensor_Z \ac Fp 
	\end{equation*}
is an \sr\ of $R$, proving the first assertion. The last assertion is now also clear.
\end{proof}

\begin{proposition}\label{P:Ftype}
Let $R$ be a local $\mathbb C$-affine domain. If $R$ has F-rational type (weakly F-regular type),  then $R$ is \genrat\ (respectively, \genreg). 
\end{proposition}
\begin{proof}
Suppose first that $R$ has F-rational type. By definition, we can find a model $(Z,R_Z)$ of $R$ such that $R_Z/\pr R_Z$ is F-rational for all maximal ideals $\pr$ of $Z$. Let $\seq\gamma p$ and $\seq Rp$ be as in \eqref{eq:seqrp} of Corollary~\ref{C:model}. Note that $\seq\gamma p(Z)$ is the residue field of $Z$ at the maximal ideal given by the kernel of $\seq\gamma p$. Hence each $R_Z\tensor_Z\seq\gamma p(Z)$ is F-rational. Since $\seq Rp$ is obtained from this by base change over the field extension $\seq\gamma p(Z) \to\ac Fp$, we get that almost all $\seq Rp$ are F-rational. Hence $R$ is \genrat\ by Proposition~\ref{P:genrat}. 

The argument for \genregity\ is the same, using Remark~\ref{R:genreg}.
\end{proof}

\begin{theorem}\label{T:Frat}
For a local $\mathbb C$-affine domain $R$, the following four statements are equivalent.
\begin{enumerate}
\item\label{l:Frat} $R$ has F-rational type.
\item\label{l:genrat} $R$ is \genrat.
\item\label{l:aicext} $R$ is $\mathcal B$-rational.
\item\label{l:ratsing} $R$ has rational singularities.
\end{enumerate}
\end{theorem}
\begin{proof}
The implication $\eqref{l:Frat}\implies \eqref{l:genrat}$ is given by Proposition~\ref{P:Ftype} and the implication  $\eqref{l:genrat}\implies\eqref{l:aicext}$  by Corollary~\ref{C:parid}. Theorem~\ref{T:ratsing} gives $\eqref{l:aicext}\implies\eqref{l:ratsing}$  and the implication $\eqref{l:ratsing}\implies\eqref{l:Frat}$ is proven by \name{Hara} in \cite{HaRat}.
\end{proof}

In particular, this proves Theorem~\ref{T:RatSing} from the introduction.  Note that \name{Smith} has already proven $\eqref{l:Frat}\implies\eqref{l:ratsing}$ in \cite{SmFrat}. Recall that we showed in \cite[Theorem 6.2]{SchRatSing} that \nsrat\ implies rational singularities. It is natural to ask whether the converse is also true.  There is another related notion which is expected to be equivalent with rational singularities, to wit, \emph{F-rationality}, that is to say, the property that some ideal generated by a system of parameters is equal to its (classical) \ch\ zero tight closure. Since \ch\ zero tight closure (more precisely, equational tight closure) is the smallest of all closure operations (see \cite[Theorem 10.4]{SchNSTC}), F-rationality is implied by $\mathcal B$-rationality. Of all implications, $\eqref{l:ratsing}\implies\eqref{l:Frat}$ is the least elementary, since \name{Hara}'s proof rests on some deep vanishing theorems.

Since rational singularities are preserved under localization, so is being $\mathcal B$-rational or being \genrat.

\begin{proposition}\label{P:Gor}
If a local $\mathbb C$-affine domain $R$ is Gorenstein and \genrat, then it is \genregglob\ whence  $\mathcal B$-regular.
\end{proposition}
\begin{proof}
As we just observed, \genratity\ is preserved under localization, so that it suffices to show that $R$ is \genreg. Let $\seq Rp$ be an   \sr\ of $R$. By \cite[Theorem 4.6]{SchNSTC},  almost all $\seq Rp$ are  Gorenstein. By Proposition~\ref{P:genrat}, almost all $\seq Rp$ are F-rational. Therefore, almost all $\seq Rp$ are F-regular, by \cite[Theorem 1.5]{HuTC}. Hence $R$ is \genreg\ by Remark~\ref{R:genreg}, whence  $\mathcal B$-regular by Corollary~\ref{C:Breg}.
\end{proof}

Recall that a \homo\ $A\to B$ is called \emph{cyclically pure}, if $IB\cap A=I$, for every ideal $I$ of $A$.

\begin{proposition}\label{P:pure}
If $R\to S$ is a cyclically pure \homo\ of local $\mathbb C$-affine domains and if $S$ is \genreg, then so is $R$. The same is true upon replacing \genreg\ by weakly $\mathcal B$-regular.
\end{proposition}
\begin{proof}
Let $\id$ be an ideal in $R$ and $z$ an element in its generic tight closure. Let $\seq Rp\to\seq Sp$ be an \sr\ of $R\to S$ (that is to say, choose \sr{s} $\seq Rp$ and $\seq Sp$ for $R$ and $S$ as well as \sr{s} for the polynomials that induce the \homo\ $R\to S$; these then induce the \homo\ $\seq Rp\to \seq Sp$, for almost all $p$; see \cite[3.2.4]{SchNSTC} for more details). Let $\seq zp$ and $\seq \id p$ be \sr{s} of $z$ and $\id$. For almost all $p$, we have  that $\seq zp$ lies in the tight closure of $\seq\id p$. By persistence (\cite[Theorem 2.3]{HuTC}), $\seq zp$ lies in the tight closure of $\seq\id p\seq Sp$, for almost all $p$, showing that $z$ lies in the generic tight closure of $\id S$. In fact, the preceding argument shows that generic tight closure is persistent (we have not yet used the purity of $R\to S$ nor even  its injectivity).  Now, by assumption, $S$ is \genreg, so that $z\in \id S$ and hence, by cyclic purity,  $z\in \id S\cap R=\id$.

To prove the last statement, observe that our assumptions imply that
	\begin{equation*}
	\id\subset \id\BCM R\cap R\subset (\id\BCM S\cap S)\cap R = \id S\cap R=\id.
	\end{equation*}
\end{proof}

\subsection{Proof of Boutot's Theorem under the additional Gorenstein hypothesis}\label{p:Bou}
Let $R\to S$ be a cyclically pure \homo\ of local $\mathbb C$-affine domains and assume $S$ is Gorenstein and has rational singularities. It follows that $S$ is $\mathcal B$-rational,  by Theorem~\ref{T:Frat}, whence weakly $\mathcal B$-regular, by Proposition~\ref{P:Gor}. Therefore, $R$ is weakly $\mathcal B$-regular by Proposition~\ref{P:pure} and hence has rational singularities by Theorem~\ref{T:Frat} again.\qed

Note that \name{Boutot}  proves the same result without the Gorenstein hypothesis.  It follows from his result that being \genrat\ (or, equivalently, being of F-rational type) descends under pure maps. However, it is not clear how to prove this from the definitions alone.

\section{\BS\ Theorems}

\subsection{Proof of Theorem~\ref{T:BS}}
Let $R$ and $I$ be as in the statement and let $z$ be an element in the integral closure of $I^{n+k}$, for some $k\in\nat$. Take \sr{s} $\seq Rp$, $\seq Ip$ and $\seq zp$ of $R$, $I$ and $z$ respectively. Since $z$ satisfies an integral equation
	\begin{equation*}
	z^n +a_1z^{n-1}+\dots+a_n=0
	\end{equation*}
with $a_i\in I^{(n+k)i}$, we have for almost all $p$ an equation
	\begin{equation*}
	(\seq zp)^n +\seq{a_1}p(\seq zp)^{n-1}+\dots+\seq{a_n}p=0
	\end{equation*}
with $\seq{a_i}p\in(\seq Ip)^{(n+k)i}$ an \sr\ of $a_i$. In other words,  $\seq zp$ lies in the integral closure of $(\seq Ip)^{n+k}$, for almost all $p$. By \cite[Theorem 7.1]{HHbigCM2},  almost all $\seq zp$ lie in $(\seq Ip)^{k+1}\seq Rp^+\cap\seq Rp$. Taking ultraproducts, we get that $z\in I^{k+1}\BCM R\cap R$, as we needed to show.\qed

In fact, the ideas in the  proof of \cite[Theorem 7.1]{HHbigCM2} can be used to carry out the argument directly in $\BCM R$. Using Theorem~\ref{T:BS}, we also get a new proof of a result of \name{Lipman} and \name{Teissier} in \cite{LT}. We need a result on powers of parameter ideals.

\begin{proposition}\label{P:powreg}
Let $R$ be a   local $\mathbb C$-affine domain with rational singularities. If  $I$ is an ideal generated by a regular sequence, then  $I^n=I^n\BCM R\cap R$, for each $n$.
\end{proposition}
\begin{proof}
Let $\mathbf x$ be a regular sequence generating $I$. We induct on $n$. If $n=1$, the assertion follows from \eqref{F:par} in Proposition~\ref{P:regseq}  since $R$ is $\mathcal B$-rational by Theorem~\ref{T:Frat}. Hence assume $n>1$ and let $a\in I^n\BCM R\cap R$. By induction, $a\in I^{n-1}$, so that $a=F(\mathbf x)$ with $F$ a homogeneous polynomial over $R$ of degree $n-1$. Since $\mathbf x$ is  a $\BCM R$-regular sequence by Theorem~\ref{T:BCMeq}, it is $\BCM R$-quasi-regular (\cite[Theorem 16.2]{Mats}). In particular $a=F(\mathbf x)\in I^n\BCM R$ implies that all coefficients of $F$ lie in $I\BCM R$, whence in $I$. Therefore, $a=F(\mathbf x)\in I^n$.
\end{proof}

\begin{remark}\label{R:mon}
More generally, we have that $J=J\BCM R\cap R$ for any ideal $J$ generated by monomials in some regular sequence $\rij xd$ such that $J$ contains a power of every $x_i$. Indeed, by \cite{EH}, any such ideal is the intersection of ideals of the form $(x_1^{t_1},\dots,x_d^{t_d})R$ for some choice of $t_i\in\nat$. Hence it suffices to prove the claim for $J$ of the latter form, and this is clear by \eqref{F:par} in Proposition~\ref{P:regseq}.
\end{remark}

\begin{theorem}[\name{Lipman-Teissier}]
If a $d$-dimensional local $\mathbb C$-affine domain $R$ has rational singularities, then for any ideal $I$ of $R$ and any $k\geq 0$, the integral closure of $I^{d+k}$ is contained in $I^{k+1}$.
\end{theorem}
\begin{proof}
 Assume first that $I$ is  generated by a system of parameters. By Theorem~\ref{T:BS},  the integral closure of $I^{d+k}$ lies in $I^{k+1}\BCM R\cap R$ and the latter ideal is just $I^{k+1}$ by Proposition~\ref{P:powreg}. Next assume that $I$ is $\maxim$-primary, where $\maxim$ denotes the maximal ideal of $R$. By \cite[Theorem 14.14]{Mats}, we can find a system of parameters $\mathbf x$ of $R$ such that $J:=\mathbf xR$ is a reduction of $I$. Since $I^{d+k}$ and $J^{d+k}$ have then the same integral closure, our previous argument shows that this integral closure lies inside $J^{k+1}$ whence inside $I^{k+1}$. Finally, let $I$ be arbitrary and put $J_n:=I+\maxim^n$. If $a$ lies in the integral closure of $I^{d+k}$, then for each $n$, it lies also in the integral closure of $J_n^{d+k}$, whence in $J_n^{k+1}$ by our previous argument. Since 
	\begin{equation*}
	J_n^{k+1}\subset  I^{k+1}+\maxim^n,
	\end{equation*}
 we get that $a$ lies in right hand side ideal for each $n$, and hence by Krull's Intersection Theorem, in $I^{k+1}$, as required.
\end{proof}

\subsection{$\mathcal B$-closure}\label{s:Bcl}
In analogy with plus closure in positive \ch\ (that is to say, the closure operation given as $I^+:=IR^+\cap R$), we define the \emph{$\mathcal B$-closure} of an ideal $I$ in a local $\mathbb C$-affine domain $R$ as the ideal $I\BCM R\cap R$. This closure operation satisfies many properties of classical tight closure, to wit: (i) a regular ring is  weakly  $\mathcal B$-regular (by Corollary~\ref{C:reg}), (ii)  colon capturing holds, in the sense that $(\rij x{i-1}R:x_iR)$ is contained in the $\mathcal B$-closure of $\rij x{i-1}R$, for every system of parameters of $R$ (by Theorem~\ref{T:BCMeq}); (iii)   \BS\ as stated in Theorem~\ref{T:BS} holds; and (iv), persistence holds (by the weak functorial property of $\BCM\cdot$). Unfortunately, in view of the more complicated way $\BCM\cdot$ and localization commute (see  Formula~\eqref{eq:loc}), it is not clear whether $\mathcal B$-closure commutes with localization (in contrast with plus closure, which is easily seen to commute with localization).

\section{Regularity and Betti numbers}

In this section, we extend the main results of \cite{SchBetti} to $\mathbb C$-affine domains. We start with proving Theorem~\ref{T:isosing} from the introduction.

\begin{theorem}\label{T:reg}
Let $(R,\maxim)$ be a local  $\mathbb C$-affine domain  with residue field $k$. If $R$ has at most an isolated singularity or has dimension at most two and if $\tor R1{\BCM R}k=0$, then $R$ is regular.
\end{theorem}
\begin{proof}
Let $(\seq Rp,\seq\maxim p)$ be  an \sr\ of $(R,\maxim)$ and let $\seq kp$ be the corresponding residue fields. It follows from  \cite[Theorems 4.5 and 4.6]{SchNSTC} that $\seq Rp$ has at most an isolated singularity or has dimension at most two, for almost all $p$. I claim that  $\tor {\seq Rp}1{\seq Rp^+}{\seq kp}=0$, for almost all $p$. Assuming the claim, we get by the Main Theorem of \cite{SchBetti} that almost all $\seq Rp$ are regular. By another application of \cite[Theorem 4.6]{SchNSTC} (see also \cite[Theorem 5.3]{SchBC}), we get that $R$ is regular, as required.

To prove the claim, we argue as follows. Write each $\seq Rp^+$ as $\pol{\seq Rp}X/\seq{\mathfrak n}p$, where $X$ is an infinite  tuple of variables and $\seq{\mathfrak n}p$ some ideal. Put $\seq Ap:=\pol{\seq Rp}X$ and let $\ul A$ and $\ul{\mathfrak n}$ be the ultraproduct of the $\seq Ap$ and the $\seq{\mathfrak n}p$ respectively. Therefore, $\BCM R=\ul A/\ul{\mathfrak n}$. The vanishing of $\tor R1{\BCM R}k$ means that $\maxim\ul A\cap\ul{\mathfrak n}=\maxim\ul{\mathfrak n}$. The vanishing of $\tor {\seq Rp}1{\seq Rp^+}{\seq kp}$ is then equivalent with the equality $\seq\maxim p\seq Ap\cap\seq{\mathfrak n}p=\seq\maxim p\seq{\mathfrak n}p$. Therefore, assume that this equality does not hold for almost all $p$, so that there exists $\seq fp$ which lies in  $\seq\maxim p\seq Ap\cap\seq{\mathfrak n}p$, but, for almost $p$ does not lie in $\seq\maxim p\seq{\mathfrak n}p$. Let $\ul f$ be the ultraproduct of the $\seq fp$. It follows from \los\ that $\ul f$ lies in $\maxim\ul A\cap\ul{\mathfrak n}$ whence in $\maxim\ul{\mathfrak n}$. Let $\maxim:=\rij ysR$ and let $\seq{y_i}p$ be an \sr\ of $y_i$, so that $\seq\maxim p=(\seq{y_1}p,\dots,\seq{y_s}p)\seq Rp$, for almost all $p$. Since $\ul f\in\maxim\ul{\mathfrak n}$, there exist $\ul{g_i}\in\ul{\mathfrak n}$, such that $\ul f=\ul{g_1}y_1+\dots+\ul{g_s}y_s$. Hence, if we choose $\seq{g_i}p\in\seq {\mathfrak n}p$ such that their ultraproduct is $\ul{g_i}$, then by \los, $\seq fp=\seq{g_1}p\seq{y_1}p+\dots+ \seq{g_s}p\seq{y_s}p$, contradicting our assumption on $\seq fp$.
\end{proof}

In general, we can prove at least the following.

\begin{corollary}
Let $R$ be a  local $\mathbb C$-affine domain  with residue field $k$. If   $\tor R1{\BCM R}k$ vanishes, then $R$ has rational singularities.
\end{corollary}
\begin{proof}
By \cite[Theorem 2.2]{SchBetti}, the vanishing of $\tor R1{\BCM R}k$ implies that $R\to\BCM R$ is cyclically pure. In particular, $\mathbf xR=\mathbf x\BCM R\cap R$ for some system of parameters $\mathbf x$  of $R$. Hence  $R$ has rational singularities by Theorem~\ref{T:ratsing}. 

Alternatively, from the proof of Theorem~\ref{T:isosing}, we get that $\tor{\seq Rp}1{\seq Rp^+}{\seq kp}$ vanishes, for almost all $p$, where $\seq Rp$ is an \sr\ of $R$ and $\seq kp$ the residue field of $\seq Rp$. By \cite[Theorem B]{SchBetti}, we get that almost all $\seq Rp$ are pseudo-rational. By \cite[Theorem 5.2]{SchRatSing}, it follows that $R$ has rational singularities.
\end{proof}

We actually showed that $R$ as above is weakly $\mathcal B$-regular.

\begin{remark}
The following are equivalent  for a local $\mathbb C$-affine domain $R$:
\begin{enumerate}
\item $R$ is regular;
\item $R\to \inverse\ulfrob{\ul R}$ is flat;
\item $\tor R1{\inverse\ulfrob{\ul R}}k=0$,
\end{enumerate}
where $k$ denotes the residue field of $R$ and where $\inverse\ulfrob{\ul R}$ is the subring of $\BCM R$ consisting of all elements whose image under $\ulfrob$ lies in $\ul R$. Indeed, let $\seq Rp$ be an \sr\ of $R$ and let $\seq kp$ the residue field of $\seq Rp$. By  Kunz's Theorem, the regularity of $\seq Rp$ is equivalent to the flatness of $\seq Rp\to \seq Rp^{1/p}$, and by the Local Flatness Criterion, this in turn is equivalent to the vanishing of $\tor {\seq Rp}1{\seq Rp^{1/p}}{\seq kp}$. Moreover, $R$ is regular \iff\ almost all $\seq Rp$ are regular (\cite[Theorem 4.6]{SchNSTC}) whereas the same argument as in the proof of Theorem~\ref{T:isosing} shows that the vanishing of  $\tor R1{\inverse\ulfrob{\ul R}}k$ is equivalent with the vanishing of almost all $\tor {\seq Rp}1{\seq Rp^{1/p}}{\seq kp}$. This proves that all assertions are equivalent.
\end{remark}

\providecommand{\bysame}{\leavevmode\hbox to3em{\hrulefill}\thinspace}

\end{document}